# SYMPLECTIC $F_2$-VECTOR SPACES AND DYCK WORDS

G. Lusztig

## Introduction

**0.1.** Let $F_2$ be the field with two elements. Let $V$ be an $F_2$-vector space of finite dimension $2n$ ($n \in \mathbf{Z}_{\geq 1}$) endowed with a nondegenerate symplectic form $<,>: V \times V \to F_2$ and with a collection of vectors $\{e_i; i \in S\}$ where $S$ is the set of vertices of a graph in which $i, j$ are joined whenever $< e_i, e_j >= 1$; it is assumed that this is a graph of affine type $A_{2n}$ (in particular $|S| = 2n+1$). We say that $\{e_i; i \in S\}$ is a circular basis of $V$ (see [L20a]). Note that $\sum_{i \in S} e_i = 0$ and any $2n$ vectors in $\{e_i; i \in S\}$ form a basis of $V$.

In this paper we show that the circular basis defines a partition of $V$ into subsets which are naturally products of sets of Dyck words (hence have cardinal equal to a product of Catalan numbers).

**0.2.** For any $I \subset S$ we denote by $\underline{I}$ the full subgraph of $S$ whose set of vertices is $I$. Let $\mathcal{I}^1$ be the set of all $I \subset S$ such that $\underline{I}$ is a graph of type $A_{2m+1}$ for some $m \in \mathbf{N}$. For $I \in \mathcal{I}^1$ let

$\delta I = \{s \in I; \underline{I - s}$ is connected or empty$\}$.

Let $I \in \mathcal{I}^1$. A map $f : I \to \mathbf{Z}_{\geq 1}$ is said to be a Dyck map if

$f(i) = 1$ for $i \in \delta I$,

$f(i) - f(i') \in \{1, -1\}$ whenever $i, i'$ are joined in $\underline{I}$.

Let $Dy(I)$ be the set of all Dyck maps $I \to \mathbf{Z}_{\geq 1}$.

For $r \in \mathbf{Z}$ we define $[r] \in F_2$ by:

$[r] = 1$ if $r = 1 \mod 4$ or $r = 2 \mod 4$,

$[r] = 0$ if $r = 3 \mod 4$ or $r = 0 \mod 4$.

Let $Dy_2(I)$ be the set of all maps $g : I \to F_2$ with the following property: there exists $f \in Dy(I)$ such that $g(i) = [f(i)]$ for all $i \in I$. The obvious surjective map $Dy(I) \to Dy_2(I)$ is denoted by $f \mapsto \bar{f}$. We show:

(a) This map is bijective.

We can assume that $S$ consists of $1, 2, 3, \dots, 2n+1$ with $(1,2)$ joined, $(2,3)$ joined, ..., $(2n, 2n+1)$ joined and $(2n_1, 1)$ joined, and that $I$ consists of $1, 2, 3, \dots, 2m+1$ with $0 \leq m < n$.

Let $f, f'$ be in $Dy(I)$ be such that $\bar{f} = \bar{f}'$. We show that $f(k) = f'(k)$ for all $k \in [1, m+1]$. We use induction on $k$.

Typeset by $\mathcal{A}\mathcal{M}\mathcal{S}$-TEX

We have $f(1)=f'(1)=1$. Assume that $k\in[2,2m]$ and that we already know that $f(k)=f'(k)$. We have $f(k+1)=f(k)+\delta, f'(k+1)=f'(k)+\delta'=f(k)+\delta'$, where $\delta=\pm1,\delta'=\pm1$ and $\bar{f}(k+1)=\bar{f}'(k+1)$ that is $[f(k+1)]=[f'(k+1]$ so that $[f(k)+\delta]=[f(k)+\delta']$. Setting $r=f(k)+\delta$ we have $[r]=[r+\delta'-\delta]$. If $\delta\ne\delta'$ we have $\delta'-\delta\in\{2,-2\}$ and using the definition we have $[r+\delta'-\delta]\ne[r]$, a contradiction. We see that $\delta=\delta'$ and $f(k+1)=f'(k+1)$ as required. This proves (a).

**0.3.** For $I,I'$ in $\mathcal{I}^1$ we write $I\spadesuit I'$ whenever $I\cap I'=\emptyset$ and $\underline{I\cup I'}$ is disconnected. Let $R$ be the set whose elements are finite unordered sequences of objects of $\mathcal{I}^1$. Let $\phi'(V)$ be the set of all $B'\in R$ such that the following holds:

$(P'_0)$ If $I\in B',I'\in B'$, then $I=I'$, or $I\spadesuit I'$.

For $B'\in\phi'(V)$ let $V^{B'}$ be the subset of $V$ consisting of all vectors of the form $\sum_{I\in B'}\sum_{i\in I}g_I(i)e_i$ for various $(g_I)\in\prod_{I\in B'}Dy_2(I)$.

**Theorem 0.4.** *We have $V=\sqcup_{B'\in\phi'(V)}V^{B'}$.*

The proof (given in §1) uses a collection $\mathcal{F}(V)$ of isotropic subspaces (introduced in [L20],[L20a]) of $V$ attached to the circular basis.

## 1. Proof of Theorem 0.4

**1.1.** We recall the definition of $\mathcal{F}(V)$ (by induction on $n\ge1$). If $n=1$, $\mathcal{F}(V)$ consists of the various subspaces of dimension $\le1$ of $V$. Assume now that $n\ge2$. For $i\in S$ let $V^i$ be the subspace of $V$ generated by the vectors $e_j$ where $j\in S$ is not joined with $i$ and by $e_j+e_i+e_{j'}$ where $j\ne j'$ are joined with $i$.

Then these vectors form a circular basis of the symplectic vector space $V^i$ with the induced $<,>$. We have $\dim V^i=2n-2$ hence $\mathcal{F}(V^i)$ is defined. A subspace of $V$ is said to be in $\mathcal{F}(V)$ if either it is 0 or if for some $i\in S$ it is of the form $E+F_2e_i$ for some $E\in\mathcal{F}(V^i)$.

For $I,I'$ in $\mathcal{I}^1$ we write $I\prec I'$ whenever $I\subset I'-\delta I'$. For $I\in\mathcal{I}^1$ let $I^{ev}$ be the set of all $s\in I$ such that $I-\{s\}=I'\sqcup I''$, with $I'\in\mathcal{I}^1$, $I''\in\mathcal{I}^1$. We have $|I^{ev}|=(|I|-1)/2$.

Let $\phi(V)$ be the set consisting of all $B\in R$ such that $(P_0),(P_1)$ below hold.

$(P_0)$ If $I\in B,I'\in B$, then $I=I'$, or $I\spadesuit I'$, or $I\prec I'$, or $I'\prec I$.

$(P_1)$ Let $I\in B$. There exist $I_1,I_2,\dots,I_k$ in $B$ such that $I^{ev}\subset I_1\cup I_2\cup\dots\cup I_k$ (disjoint union), $I_1\prec I$, $I_2\prec I,\dots,I_k\prec I$.

For $I\in\mathcal{I}^1$ let $e_I=\sum_{i\in I}e_i\in V$. For $B\in\phi(V)$ let $L_B$ the subspace of $V$ spanned by $\{e_I;I\in B\}$. The following is proved in [L24].

(a) *For any $B\in\phi(V)$, the set $\{e_I;I\in B\}$ is a basis of $L_B$. The map $B\mapsto L_B$ is a bijection $\phi(V)\to\mathcal{F}(V)$ whose inverse maps any $L\in\mathcal{F}(V)$ to the set $\{I\in\mathcal{I}^1;e_I\in L\}$ (which is in $\phi(V)$).*

**1.2.** Let $B\in\phi(V)$. Let ${}^0B$ be the set of all $I\in B$ such that $I$ is not properly contained in any $I'\in B-\{I\}$. We have ${}^0B\in\phi'(V)$.

We define

$$\tilde{\phi}(V) = \{(B', (f_I)_{I \in B'}); B' \in \phi'(V), f_I \in Dy(I)\}.$$

If $B \in \phi(V)$ we define $f^B : S \to \mathbf{Z}_{\ge 1}$ by $f^B(i) = |I \in B; i \in I|$. One can verify (using the definitions) that

(a) For any $I \in {}^0B$ we have $f^B|_I \in Dy(I)$

and that

(b) *the map* $\phi(V) \to \tilde{\phi}(V)$ *given by associating to* $B \in \phi(V)$ *the subset* ${}^0B$ *and the collection of Dyck maps* $f^B|_I \in Dy(I)$ *(for various* $I \in {}^0B$*) is a bijection. The inverse map* $\tilde{\phi}(V) \to \phi(V)$ *associates to* $(B', (f_I)_{i \in B'})$ *the set of connected components of* $\underline{\{i \in I; f_I(i) = k\}}$ *for various* $I \in B'$ *and various* $k \in \mathbf{Z}_{\ge 1}$.

**1.3. Examples.** In this subsection we assume that $B' = \{I\} \in \phi'(V)$ where $I \in \mathcal{I}^1$ is such that $|I| \in \{3, 5, 7, 9\}$. For such $I$ and for a fixed $f_I \in Dy(I)$ we describe a matrix with some entries marked by $X$ and with other entries unmarked. The columns of the matrix are indexed by the elements of $I$ so that two adjacent columns correspond to two elements of $I$ which are joined. The column indexed by $i \in I$ has $f_I(i)$ marked entries. The elements of the corresponding $B$ can be viewed as the connected components of the set of marked entries in the various rows of the matrix. In each case we write the sequence $f$ of numbers $f_I(i)$ (with $i \in I$).

$$\begin{pmatrix} X & X & X \\ & X & \end{pmatrix} ...f = 121$$

$$\begin{pmatrix} X & X & X & X & X \\ & X & X & X & \\ & & X & & \end{pmatrix} ...f = 12321$$

$$\begin{pmatrix} X & X & X & X & X \\ & X & & X & \end{pmatrix} ...f = 12121$$

$$\begin{pmatrix} X & X & X & X & X & X & X \\ & X & X & X & X & X & \\ & & X & X & X & & \\ & & & X & & & \end{pmatrix} ...f = 1234321$$

$$\begin{pmatrix} X & X & X & X & X & X & X \\ & X & X & X & X & X & \\ & & X & & X & & \end{pmatrix} ...f = 1232321$$

$$\begin{pmatrix} X & X & X & X & X & X & X \\ & X & X & X & & X & \\ & & X & & & & \end{pmatrix} ...f = 1232121$$

$$\begin{pmatrix} X & X & X & X & X & X & X \\ & X & & X & X & X & \\ & & & & X & & \end{pmatrix} ...f = 1212321$$

$$\begin{pmatrix} X & X & X & X & X & X & X \\ & X & & X & & X & \end{pmatrix} ...f = 1212121$$

$$\begin{pmatrix} X & X & X & X & X & X & X & X & X \\ & X & X & X & X & X & X & X & \\ & & X & X & X & X & X & & \\ & & & X & X & X & & & \\ & & & & X & & & & \end{pmatrix} ...f = 123454321$$

$$\begin{pmatrix} X & X & X & X & X & X & X & X & X \\ & X & X & X & X & X & X & X & \\ & & X & X & X & X & X & & \\ & & & X & & X & & & \end{pmatrix} ...f = 123434321$$

$$\begin{pmatrix} X & X & X & X & X & X & X & X & X \\ & X & X & X & X & X & X & X & \\ & & X & X & X & & X & & \\ & & & X & & & & & \end{pmatrix} ...f = 123432321$$

$$\begin{pmatrix} X & X & X & X & X & X & X & X & X \\ & X & X & X & X & X & X & X & \\ & & X & & X & X & X & & \\ & & & & & X & & & \end{pmatrix} ...f = 123234321$$

$$\begin{pmatrix} X & X & X & X & X & X & X & X & X \\ & X & X & X & X & X & X & X & \\ & & X & & X & & X & & \end{pmatrix} ...f = 123232321$$

$$\begin{pmatrix} X & X & X & X & X & X & X & X & X \\ & X & X & X & X & X & & X & \\ & & X & X & X & & & & \\ & & & X & & & & & \end{pmatrix} ...f = 123432121$$

$$\begin{pmatrix} X & X & X & X & X & X & X & X & X \\ & X & X & X & X & X & & X & \\ & & X & & X & & & & \end{pmatrix} ...f = 123232121$$

$$\begin{pmatrix} X & X & X & X & X & X & X & X & X \\ & X & & X & X & X & X & X & \\ & & & & X & X & X & & \\ & & & & & X & & & \end{pmatrix} ...f = 121234321$$

$$\begin{pmatrix} X & X & X & X & X & X & X & X & X \\ & X & & X & X & X & X & X & \\ & & & & X & & X & & \end{pmatrix} ...f = 121232321$$

$$\begin{pmatrix} X & X & X & X & X & X & X & X & X \\ & X & X & X & & X & X & X & \\ & & X & & & & X & & \end{pmatrix} ...f = 123212321$$

$$\begin{pmatrix} X & X & X & X & X & X & X & X & X \\ & X & X & X & & X & & X & \\ & & X & & & & & & \end{pmatrix} ...f = 123212121$$

$$\begin{pmatrix} X & X & X & X & X & X & X & X & X \\ & X & & X & X & X & & X & \\ & & & & X & & & & \end{pmatrix} ...f = 121232121$$

$$\begin{pmatrix} X & X & X & X & X & X & X & X & X \\ & X & & X & & X & X & X & \\ & & & & & & X & & \end{pmatrix} ...f = 121212321$$

$$\begin{pmatrix} X & X & X & X & X & X & X & X & X \\ & X & & X & & X & & X & \end{pmatrix} ...f = 121212121$$

**1.4.** According to [L24],

(a) the map $B \mapsto \sum_{i \in S} [f^B(i)] e_i$ is a bijection $\mathcal{F}(V) \to V$.

Using this and 1.2(b) we deduce that

(b) the map $(B', (f_I)_{I \in B'}) \mapsto \sum_{I \in B'} \sum_{i \in I} [f_I(i)] e_i$ is a bijection $\tilde{\mathcal{F}}(V) \to V$.

Now Theorem 0.4 follows immediately from (b).

## 2. Dyck words

**2.1.** Recall that a Dyck word of length $2m$ is string containing $m$ letters $a$ and $m$ letters $b$ in such a way that for any $k \geq 1$ the first $k$ letters of the string involve at least as many $a$ as $b$. Here are the Dyck words of length $2, 4, 6$.

$ab$

$aabb, abab$

$aaabbb, aababb, aabbab, abaabb, ababab.$

It is known that the number of Dyck words of length $2m$ is the Catalan number $(m+1)^{-1}\binom{2m}{m}$.

Let $I \in \mathcal{I}^1$. We write $I = \{1, 2, \dots, 2m+1\}$ as in the proof of 0.2(a). For any $f \in Dy(I)$ we define a string of letters $a, b$ as follows. For $k \in [1, 2m]$ the $k$-th letter of the string is $a$ if $f(k+1) = f(k) + 1$ and is $b$ if $f(k+1) = f(k) - 1$. We thus obtain a Dyck word of legth $2m$. This gives a bijection between $Dy(I)$ and the set of Dyck words of legth $2m$. It follows that $|Dy(I)| = |Dy_2(I)|$ is a Catalan number.

For $I$ as above, the set $Dy(I)$ is
121 if $m=1$
$12321, 12121$ if $m=2$
$1234321, 1232321, 1232121, 1212321, 1212121$ if $m=3$

$$123454321, 123434321, 123432321, 123234321, 123432121, 123232321, 121234321,$$
$$123232121, 121232321, 123212321, 123212121, 121232121, 121212321, 121212121$$

if $m=4$.

The corresponding elements of $Dy_2(I)$ are
111
$11011, 11111$
$1100011, 1101011, 1101111, 1111011, 1111111$

$$110010011, 110000011, 110001011, 110100011, 110001111, 110101011, 111100011,$$
$$110101111, 111101011, 110111011, 110111111, 111101111, 111111011, 111111111$$

Department of Mathematics, M.I.T., Cambridge, MA 02139